\documentclass{amsart}

\usepackage{amssymb}  
\usepackage{amsmath}  
\usepackage{latexsym} 

\usepackage{hyperref}

\newcommand{\be}{\begin{equation}}
\newcommand{\ee}{\end{equation}}
\newcommand{\bea}{\begin{eqnarray}}
\newcommand{\eea}{\end{eqnarray}}

\begin{document}

\title{a Characterization of Koiso's typed Solitons}

\author{yang bo}
\address{School of Mathematics Sciences, Peking
University, Beijing, 100871, People's Republic of China.}
\email{{yangbopku@gmail.com}}

\begin{abstract}
By extending Koiso's examples to the non-compact case, we construct
complete gradient K\"{a}hler-Ricci solitons of various types on
certain holomorphic line bundles over compact K\"{a}hler-Einstein
manifolds. Moreover, a uniformization result on steady gradient
K\"{a}hler-Ricci solitons with non-negative Ricci curvature is
obtained under additional assumptions.
\end{abstract}

\maketitle

\section{Introduction}

A K\"{a}hler metric $g$ is called a K\"{a}hler-Ricci soliton on a
complex manifold $M$ if there is a holomorphic vector field $V$
and a real number $\rho$ such that the equality $Ric+\rho
g-L_{V}g=0$ holds on $M$. It is called steady when $\rho=0$,
expanding when $\rho>0$ and shrinking when $\rho<0$. In addition,
when $V$ is the gradient vector field of a real function $f$ on
$M$ we call it a gradient K\"{a}hler-Ricci soliton.

Like Ricci solitons in the real case, K\"{a}hler-Ricci solitons
naturally arise as one takes limits of dilations of singularities in
K\"{a}hler-Ricci flow (cf. \cite{H2} and \cite{C2}). Due to their
importance in the study of K\"{a}hler-Ricci flow, it is interesting
to learn more examples of K\"{a}hler-Ricci solitons. One typical
method of constructing K\"{a}hler-Ricci solitons is to impose
certain symmetry conditions and reduce the soliton equation to ODEs
which are more tractable (for instance, in \cite{K}, \cite{C1} for
the compact case and \cite{H1},\cite{C1},\cite{C2},\cite{F-I-K} for
the non-compact case.). In the other direction, Wang and Zhu
\cite{W-Z} constructed K\"{a}hler-Ricci solitons on toric K\"{a}hler
manifolds with positive first Chern class by solving equations of
$\mathrm{Monge\!\!-\!\!Amp\grave{e}re}$ type on toric Fano
manifolds.

To the author's knowledge, all known examples in the non-compact
case rely on the $\mathrm  {U(n)}$-symmetry assumption of the
soliton metric. However, it is easy to note that Koiso's methods in
\cite{K}, which allows less restriction on the symmetry assumption,
may work similarly in the non-compact case. Since there seems to be
no literature containing this explicitly, in this note we carry out
this construction in details. We prove:

\vskip  2 mm

\noindent \textbf{Theorem 1.1}. \textit{Let M be a compact
$(n-1)$-dim K\"{a}hler-Einstein manifold satisfying $Ric(g_0)=g_0$
and L be a holomorphic line bundle over M. Assume there exists a
Hermitian metric on L such that the eigenvalues $\lambda_{i}$ of the
$Ricci$ form of $L$ with respect to $g_{0}$ are constant on $M$.
Then We have:}

\textit{(1) If $-1<\lambda_{i}<0$ for $1\leqslant i \leqslant n-1$,
there exists a complete shrinking gradient K\"{a}hler-Ricci soliton
on L.}

\textit{(2) If $\lambda_{i}<-1$ for $1\leqslant i \leqslant n-1$,
there exists a family of complete expanding gradient
K\"{a}hler-Ricci solitons on L.}

\textit{(3) If $L$ is the canonical line bundle, there exists a
family of complete steady gradient K\"{a}hler-Ricci solitons on L.}

\vskip  2 mm

Rotationally symmetric K\"{a}hler-Ricci solitons on holomorphic line
bundles over $\mathrm{CP^{n-1}}$ were constructed in \cite{C1} and
\cite{F-I-K}. It can be checked that Theorem 1.1 includes their
examples as a sepecial case. Moreover, in the steady and expanding
case, Theorem 1.1 provides a family of K\"{a}hler-Ricci solitons
over certain holomorphic line bundles of compact K\"{a}hler-Einstein
manifolds. Those line bundles can be chosen to be canonical line
bundles or their tensor products.

One can check that the steady gradient K\"{a}hler-Ricci solitons on
canonical line bundles over K\"{a}hler-Einstein manifolds in Theorem
1.1 have nonnegative Ricci curvature on the zero section and
positive Ricci curvature away from it. This motivates us to prove
the following uniformization result as a characterization of Koiso's
typed soliton in the steady case.

\vskip  2 mm

\noindent \textbf{Theorem 1.2}. \textit{Let M be a non-compact
steady gradient K\"{a}hler-Ricci soliton with non-negative Ricci
curvature. Assume its scalar curvature attains a positive maximum
along a compact complex submanifold K with codimension 1 and the
Ricci curvature is positive away from K. Then M is biholomorphic to
a holomorphic line bundle over K.}

\vskip  2 mm

It was shown in \cite{B} and \cite{C-T} independently that a
non-compact complex manifold is biholomorphic to $\mathrm C^n$ if it
admits a steady gradient K\"{a}hler-Ricci soliton metric with
positive Ricci curvature and its scalar curvature attaining a
maximum at some point. In \cite{B} Bryant proved a crucial lemma on
the existence of a nice local coordinate near singular points of the
gradient holomorphic vector field associated to the soliton and
obtained his uniformization result. Here we apply Bryant's lemma in
\cite{B} in our case to prove Theorem 1.2.

It is also interesting to note that some results concerning the
geometric classification of shrinking Ricci solitons under curvature
assumptions was shown in \cite{N} and the recent preprints
\cite{N-W},\cite{P-W} and \cite{Na}.

\vskip  2 mm

\noindent \textbf{Acknowledgments}: The author thanks Professor
Xiaohua Zhu for his guidance and encouragement, as well as Professor
Albert Chau for helpful discussions during the Summer School in
Geometric Analysis held in USTC in July 2007.

\section{Example of K\"{a}hler-Ricci solitons}

\part*{1 \ \ Background and calculations}

In this part, we quickly review some facts on Koiso's construction
of a K\"{a}hler metric on a $\mathrm {C^{\ast}-}$ bundle over a
compact K$\ddot{a}$hler manifold.

Given a holomorphic line bundle $L \rightarrow M$ on a complex
manifold M where $\pi$ is the natural projection, assume its local
trivialization is given by $\varphi_{\alpha}: U_{\alpha} \times
\mathrm{C} \rightarrow \pi^{-1}(U_{\alpha})$. We define $\mathrm
{C^{\ast}-}$ action on $L^{\ast}=L\setminus{\{0-section}\}$ by
$\varphi_{\alpha}(p,\lambda) \circ g= \varphi_{\alpha}(p,\lambda g)
\ \ \ (\forall g \in \mathrm {C^{\ast}}=C \setminus\{0\}, p\in M,
\lambda \in \mathrm C)$. It can be checked that this definition is
independent of the choice of local trivialization and this $\mathrm
{C^{\ast}=R^{+}\times S^{1}}$ action is free. We denote the two
holomorphic vector fields generated by $\mathrm{R^{+}}$ and
$\mathrm{S^{1}}$ action $H$ and $S$.

Let $L$ be a Hermitian Line bundle over a compact K\"{a}hler
manifold $M$. Denote $\tilde{J}$ the complex structure on $L$ and
$\rho_{L}$ the $Ricci$ form of $L$. Assume $t$ is a function on $L$
depending only on the norm and increasing on the norm. We consider a
hermitian metric on $L^{\ast}$ of the form
$$\tilde{g}=\pi^{\ast}g_{t}+dt^{2}+(dt\circ \tilde{J})^{2},$$
where $g_t$ is a family of $Riemannian$ metrics on $M$. Denote
$u(t)^{2}=\tilde{g}(H,H)$. It can checked that $u$ depends only on
t.

The following facts can be found in \cite{K-S}.

\vskip 2mm

\noindent \textbf{Fact 2.1}. \textit{$\tilde{g}$ is K\"{a}hler on
$L^{\ast}$ if and only if each $g_{t}$ is K\"{a}hler on $M$ and
$g_{t}=g_{0}-UB $ where $B(JX,Y)=\rho_{L}(X,Y)\ \ U=\int_{0}^{t}
u(t)dt$.}

\vskip 2mm

We further assume the eigenvalues of B with respect to $g_{0}$ are
constant on M. Let $z^{1}\cdots z^{n}$ be local coordinates of $M$
and denote $z^{0}\cdots z^{n}$ be local coordinates of $L_{0}$ such
that $\frac{\partial}{\partial z^{0}}=H-\sqrt{-1}S$.

\vskip 2mm

\noindent \textbf{Fact 2.2}. \textit{$\tilde{g}_{00}=2u^2,\
 \tilde{g}_{\alpha 0}=2u\partial_{\alpha} t,\ \tilde{g}_{\alpha
\bar{\beta}}=g_{t \alpha \bar{\beta}}+2\partial_{\alpha}
t\partial_{\bar{\beta}} t.$} \textit{Define $p=det(g_{0}^{-1}\cdot
g_{t})$, then $det(\tilde{g})=2u^2 \cdot p \cdot det(g_{0})$.}

\vskip 2mm

\noindent \textbf{Fact 2.3}. \textit{If we assume that
$\partial_{\alpha} t=\partial_{\bar{\alpha}} t=0\ \ (1 \leqslant
\alpha \leqslant n-1)$ on a fiber. If a function f is defined on
$L^{\ast}$ and only depending on $t$, then
$\partial_{0}\partial_{\bar{0}}f=u\frac{d}{dt}(u\frac{df}{dt}),\
\partial_{\alpha}\partial_{\bar{0}}f=0,\
\partial_{\alpha}\partial_{\bar{\beta}}f=-\frac{1}{2}u\frac{df}{dt} B_{\alpha
\bar{\beta}}$.}

\vskip 2mm

\noindent \textbf{Fact 2.4}. \textit{Under the same assumption of
Fact 2.3. the Ricci curvature of $\tilde{g}$ becomes:
$\tilde{R}_{00}=-u\cdot \frac{d}{dt}(u
\cdot\frac{d}{dt}(log(u^{2}p))),\ \tilde{R}_{\alpha0}=0, \ \
\tilde{R}_{\alpha \bar{\beta}}=R_{0\alpha
\bar{\beta}}+\frac{1}{2}u\cdot \frac{d}{dt}(log(u^{2}p))\cdot
B_{\alpha \beta}$.}

\vskip 2mm

Reset $\phi(U)=u(t)^2,\ Q(U)=p$ and one can compute the following:

\vskip 2mm

\noindent \textbf{Fact 2.5}. \textit{Set $V=-\frac{E}{2}H$, under
the same assumption of Fact 2.3, one can compute:}

\textit{$\tilde{R}_{00}-\tilde{g}_{00}-\textit{L}_{V}{\tilde(g)}_{00}=-H(P\circ
\phi-E\phi),\ $
$\tilde{R}_{\alpha\bar{0}}=\tilde{g}_{\alpha\bar{0}}=\textit{L}_{V}{\tilde(g)}_{\alpha\bar{0}}=0,\
$ $\tilde{R}_{\alpha \bar{\beta}}-\tilde{g}_{\alpha
\bar{\beta}}-\textit{L}_{V}{\tilde(g)}_{\alpha
\bar{\beta}}=\frac{1}{2}\ [P\circ \phi-E \phi]B_{\alpha
\bar{\beta}}+(R_{0\alpha \bar{\beta}}-g_{0\alpha \bar{\beta}}),$}
\textit{where $P \circ
\phi=\frac{d\phi}{dU}+\frac{\phi}{Q}\frac{dQ}{dU}+2U$.}

\vskip 2mm

If the initial metric $g_{0}$ on M is K\"{a}hler-Einstein satisfying
$Ric(g_0)=g_{0}$, then from the above one can reduce the shrinking
soliton equation $\tilde{Ric}-\tilde{g}-\textit{L}_{V}{\tilde(g)}=0$
to a one order ODE and get the formal solution
$$\phi(U)=-\frac{2e^{EU}}{Q(U)}\int_{Umin}^{U}xe^{-Ex}Q(x)\ dx,$$
where we denote $[Umin,Umax]$ and $[tmin,tmax]$ to be the range of
of the function $U$ and $t$ respectively. We assume that $Umin\neq
-\infty$ and $tmin \neq -\infty$.

We state the following lemma concerning the condition on grwoth of
$\phi(U)$ in order to get a well-defined K\"{a}hler metric on
${L^{\ast}}$.

\vskip  1 mm

\noindent \textbf{Lemma 2.1}. \textit{If $\phi(U)>0$ and $g_0-UB$
remains positive on $(Umin,Umax)$, in addition, $\int_{Umin}^{U}
\frac{dU}{\phi(U)}=+\infty$, $\int_{U}^{Umax}
\frac{dU}{\phi(U)}=+\infty$ and $\int_{Umin}^{U}
\frac{dU}{\sqrt{\phi(U)}}$ is finite for all $U \in (Umin,Umax)$,
then we can get a unique expression of $t$ w.r.t. the Hermitian
metric $r$ on $L$ with given initial value $tmin$, this results a
K\"{a}hler metric on $L^{\ast}$ which satisfies soliton equation on
$L^{\ast}$. } \vskip  2 mm

\begin{proof}
According to definition we know
$$\frac{dU}{\sqrt{\phi(U)}}=dt$$
$$\sqrt{\phi(U)}=u(t)=r\frac{dt}{dr}$$

this implies:
$$\int_{Umin}^{U} \frac{dU}{\phi(U)}=\int_{0}^{r} \frac{dr}{r}$$
$$\int_{Umin}^{U} \frac{dU}{\sqrt{\phi(U)}}=\int_{tmin}^{t} dt$$

Clearly when the assumption in the lemma holds, we can solve the
ODEs to get the expression of $t$ in terms of $r$ from the above
formula and $r$ varies from $0$ to $+\infty$.
\end{proof}

We remark that by solving $\frac{df}{dt}=u$ one can also make the
holomorphic vector $V$ be given by a gradient vector field of a
real-valued function $f$ which depends only on t.

\part*{2\ \ \ The shrinking case}

For convenience we introduce the following assumption.

\vskip  2 mm

\noindent \textbf{Assumption 2.1}. \textit{$\pi: L\rightarrow M$ is
a Hermitian holomorphic line bundle over a compact $(n-1)$-dim
K\"{a}hler-Einstein manifold $M$ with $Ric(g_0)=g_0$, where the
eigenvalues of the $Ricci$ form of $L$ with respect to $g_{0}$ are
constant on $M$ and satisfying $-1<\lambda_{i}<0$ for $1\leqslant i
\leqslant n-1$.}

\vskip  2 mm

Solving $\tilde{Ric}-\tilde{g}-\textit{L}_{V}{\tilde{g}}=0$, we
get the formal expression:
$$\phi(U)=-\frac{2e^{EU}}{Q(U)}\int_{Umin}^{U}xe^{-Ex}Q(x)\ dx      $$
on $L^{\ast}$. According to Assumption 2.1 the eigenvalue of $B$
with respect to $g_0$ are $-1<\lambda_{1} \cdots \lambda_{n-1}<0$,
then $Q(U)=\prod _{i=1}^{n-1} (1-U \lambda_{i})$, now the above
expression of $\phi(U)$ can be computed explicitly:
$$\phi(U)=\frac{2\eta(U,E)}{Q(U)}-\frac{2 e^{E(U-Umin)}}{Q(U)} \eta(Umin,E),$$
where $\eta(U,E)$ is a degree-n polynomial with respect to U with
the principal term $\frac{1}{E} U^n$.

Our goal is to discuss the possibility to make it complete. If it is
incomplete along the zero section we only consider adding $M$ to
complete it, and if it is still incomplete along infinity, we also
consider adding $M$ along infinity hence compactify L to a
projective bundle (see \cite{F-I-K} for other ways to complete the
metric.).

If we require $Umin=-1$, following the computation on the first
Chern class on the zero section in Example 2.1 in \cite{F-I-K}, we
can check that this is a necessary condition to extend the metric to
zero section. Calculating the metric $\tilde{g}$ in the local
coordinates one can also check that this suffices to extend the
metric non-degenerate along the zero section. A similar analysis
shows that $Umax=1$ will suffice to extend metric $\tilde{g}$
non-degenerate at infinity if we want to compactify it to a
projective bundle.

In order to get a complete non-compact metric we have to learn more
about the behavior of $\tilde{g}$ along infinity. One can also check
that if the metric is in the form of
$\tilde{g}=\pi^{\ast}g_{t}+dt^{2}+(dt\circ \tilde{J})^{2}$ the
geodesic starting in the fibre direction moves along the holomorphic
vector field $\nabla f$ away from zero section. In order to make the
soliton metric complete along infinity, one only need to check
whether those geodesic tends to infinity. This further means the
growth of those geodesics is reflected in the growth of function
$t$. (i.e. $t-tmin$ measures the length of those geodesics).

It turns out that value of E determines the behavior of this metric
along infinity. First we introduce two values of E which are
critical in our analysis. Define $E_0$ be the solution of an
algebraic equation $\eta(-1,E)=0$ and $E_1$ to be the solution to
$$\phi(1)=\frac{2\eta(1,E)}{Q(1)}-\frac{2 e^{2E}}{Q(1)} \eta(-1,E)=0$$

We have the following lemma to ensure the existence and uniqueness
of $E_1$ and $E_0$:

\vskip  2 mm \noindent \textbf{Lemma 2.2}. \textit{For any
$-1<\lambda_{1} \cdots \lambda_{n-1}<0$, $E_0$ and $E_1$ exist
uniquely. And $0<E_1<E_0<+\infty$.} \vskip  2 mm

\begin{proof}
This can be shown by writing those polynomial equations explicitly
and analyzing the signs of coefficients carefully.
\end{proof}

We now begin to analyze how E affects the asymptotic behavior of the
metric in details.

\emph{(Case 1)}: When $E=E_0$, $\phi(U)$ satisfies the assumptions
in Lemma 2.1, and one can check that when $r$ changes from $0$ to
$+\infty$, $t$ changes from $tmin$ to $+\infty$, $U$ from $Umin=-1$
to $+\infty$ and $\phi(U)$ from $0$ to $+\infty$. The resulting
metric is complete on the total space of $L$.

\emph{(Case 2)}: When $E>E_{0}$, $\phi(U)$ does not satisfy the
assumptions in Lemma 2.1. When $U$ changes from $Umin=-1$ to
$+\infty$, we have $t$ changes from $tmin$ to $tmax$ with finite
value and $\phi(U)$ from $0$ to $+\infty$, however $r$ changes
from $0$ to a finite value. This results a metric which is not
well defined on $L^{\ast}$.

\emph{(Case 3)}: When $E=E_{1}$, when $r$ changes from $0$ to
$+\infty$, $t$ changes from $tmin$ to $tmax$ with finite value, $U$
from $Umin=-1$ to $Umax=1$ and $\phi(U)$ from $0$ to $0$. One can
complete the metric by compactifying L to a projective bundle and
this results a compact gradient K\"{a}her-Ricci soliton on P(L).
This has been obtained in \cite{K} and \cite{C1}. It is interesting
to note that $E_1$ is related to the holomorphic invariant defined
in \cite{T-Z}.

\emph{(Case 4)}: For all other E, $\phi(U)$ satisfies the
assumptions in Lemma 2.1. When $r$ changes from $0$ to $+\infty$,
$t$ changes from $tmin$ to $tmax$ with finite value, $U$ from
$Umin=-1$ to a finite $Umax \neq 1$ and $\phi(U)$ from $0$ to $0$.
The resulting metric is incomplete along infinity. Moreover, it can
not be completed by adding a M along infinity since $Umin \neq 1$.

To sum up the above discussion, the shrinking case in Theorem 1.1 is
restated as follows.

\vskip 2mm

\noindent \textbf{Theorem 2.1}. \textit{Under the Assumption 2.1,
then there exists a complete shrinking Kahler-Ricci soliton on L and
the projectified line bundle P(L) respectively satisfying
$\tilde{Ric}-\tilde{g}-\textit{L}_{V}{\tilde{g}}=0$ such that
$\tilde{g}$ is in the form
$\tilde{g}=\pi^{\ast}g_{t}+dt^{2}+(dt\circ \tilde{J})^{2}$, here $V$
can be uniquely determined by the natural holomorphic $\mathrm
{R^+}$ action.}

\vskip 2mm

\part*{3\ \ \ The expanding case and steady case}

\vskip 1mm

\noindent \textbf{Assumption 2.2}. \textit{The only difference with
Assumption 2.1 is the eigenvalues of the $Ricci$ form of $L$ with
respect to $g_{0}$ are constant on $M$ and satisfying
$\lambda_{i}<-1$ for $1\leqslant i \leqslant n-1$.}

\vskip 2mm

Solving $\tilde{Ric}+\tilde{g}-\textit{L}_{V}{\tilde{g}}=0$, we
get the formal expression:
$$\phi(U)=\frac{2e^{EU}}{Q(U)}\int_{Umin}^{U}xe^{-Ex}Q(x)\ dx$$
on $L^{\ast}$.

The analysis of the asymptotic behavior of the metric is similar
to shrinking case, we only list the result:

(1) if we require $Umin=1$, we can show that this will suffice to
extend the metric to the zero section non-degenerate.

(2) Compared with the shrinking case, the restriction on E is
moderate in the expanding case.

\emph{(Case 1)}: for any $E<0$, it can be show that the resulting
metric is complete along infinity.

\emph{(Case 2)}: for any $E>0$, when $U$ from $Umin=-1$ to
$+\infty$, we have $t$ changes from $tmin$ to $tmax$ with finite
value and $\phi(U)$ from $0$ to $+\infty$, however $r$ changes
from $0$ to a finite value. This results a metric which is not
well defined on $L^{\ast}$.

\vskip  2 mm

We now turn to steady case.

\vskip  2 mm

\noindent \textbf{Assumption 2.3}. \textit{$L$ is the canonical line
bundle a compact $(n-1)$-dim $K\ddot{a}hler-Einstein$ manifold $M$
with $Ric(g_0)=g_0$.}

\vskip  2 mm

Solving $\tilde{Ric}-\textit{L}_{V}{\tilde{g}}=0$ one can get the
formal expression is:
$$\phi(U)=\frac{2e^{EU}}{Q(U)}\int_{Umin}^{U} e^{-Ex}Q(x)\ dx$$
on $L^{\ast}$.

The analysis of the asymptotic behavior of the metric shows that
the following holds:

(1) we only need to require $Umin>-1$. one can show that this will
suffice to extend the metric to the zero section non-degenerate.

(2) the result on $E$ is similar to the expanding case.

\emph{(Case 1)}: for all $E<0$, it can be show that the resulting
metric is complete along infinity. What is interesting is that in
this case the resulting metric has positive Ricci curvature away
from the zero section and nonnegative on the zero section.
Unfortunately this metric can not have nonnegative bisectional
curvature everywhere, the author thanks Prof. \hskip -1mm Albert
Chau and Prof. \hskip -1mm Fangyang Zheng for providing this
information. In fact, since the zero section is totally geodesic one
can verify that at any point on the zero section the holomorphic
bisectional curvature of the plane by the fiber direction and the
tangent direction is always negative using the curvature formula in
submanifold geometry.

\emph{(Case 2)}: if $E>0$, then this results a metric which is not
well defined even on $L^{\ast}$.

We now restate Theorem 1 in steady and expanding case as follows:

\vskip  2 mm

\noindent \textbf{Theorem 2.2}. \textit{Under Assumption 2.2 (or
Assumption 2.3), then we can find a family of complete expanding
Kahler-Ricci solitons (or complete steady Kahler-Ricci solitons) on
L where $\tilde{g}$ is in the form of
$\tilde{g}=\pi^{\ast}g_{t}+dt^{2}+(dt\circ \tilde{J})^{2}$.
Moreover, in the steady case the soliton metrics have positive Ricci
curvature away from the zero section.}

\vskip  2 mm

\section{a uniformization theorem}

In this section, we want to prove Theorem 1.2. First we state the
following lemma due to Bryant in \cite{B}.

\vskip  1 mm

\noindent \textbf{Lemma 3.1}. \textit{Let Z be the holomorphic
vector field associated to a gradient K\"{a}hler-Ricci soliton and p
be one of singular point of Z, there exists a p-centered holomorphic
coordinates $w^{1}\cdots w^{n}$ on a neighborhood $U_p$ on which
$Z=Ric(p,\frac{\partial}{\partial w^{i}})w^{i}
\frac{\partial}{\partial w^{i}}$ and $\frac{\partial}{\partial
w^{1}} \cdots \frac{\partial}{\partial w^{n}}$ are orthogonal at p.}

\vskip  2 mm

\noindent \textit{Proof of Theorem 1.2.} First we show that $\nabla
f$ vanishes on K, hence K is a totally geodesic complex submanifold.

This can be proved by the similar method in \cite{H2}. One only need
to rule out the possibility that the integral curve generated by
$\nabla f$ from any point p on K may always stay in K if $\nabla f$
does not vanish at p. This is impossible because:
$$\frac{d f(\phi_t(p))}{dt}=-{|\nabla f|}^2,$$
$$R+{|\nabla f|}^2=const,$$
where $\phi_t$ is the holomorphic flow generated by $-\nabla f$.
Then $f$ decays to infinity on the compact set K if $\phi_t(p)$
always stay in K. This contradiction shows this integral curve
must leave K after a finite time, then the remaining argument is
similar in \cite{H2}.

We further show that for any point in M, it converges to K under the
holomorphic flow $\phi_t$ generated by $-\nabla f$. Since $\nabla f$
vanishes only on K since $f$ is strictly convex outside K and weakly
convex on M. we know that $\nabla f$ vanishes only on K.

Since
$$\frac{dR(\phi_t(x))}{dt}=R_{i\overline{j}}\nabla_{i}f
\nabla_{\overline{j}}f,$$ and $dis(\phi_t(x),K)$ is non-increasing
when t increases. we conclude that $dist(\phi_t(x),K)$ converges to
zero as t goes to infinity.

We applied Lemma 3.1 under the assumption of Theorem 3.1. For an
arbitrary point p on K, the Ricci curvature of M vanishes when
restricted to $T_{p}K$. From the original proof of Lemma 3.1 in
\cite{B} we can also pick the coordinates $w^{1}\cdots w^{n}$ on a
neighborhood $U_p$ such that $\frac{\partial}{\partial w^{1}} \cdots
\frac{\partial}{\partial w^{n-1}}$ at p lie in $T_p K$. Now we find
a local holomorphic coordinates $w$ on $U_p$ which satisfies
$$w^{i}(exp_{t\nabla f}(q))=w^{i}(q)$$ when $1 \leqslant i \leqslant
n-1$ and $$w^{n}(exp_{t\nabla f}(q))=\exp(ht)w^{n}(q)$$ where $q$ is
any point in $U_{p}$ and $h=Ric(p,\frac{\partial}{\partial w^{n}})$
is a real constant.

Define $$W_{p}=\{x \mid dist(\phi_t(x),p)\rightarrow 0 \ \ \ \ \ as
\ \ \ t \rightarrow +\infty\},$$ then M can be written as
$\bigcup_{p \in K}{W_p}$. We now show that the above coordinates
$w^{1}\cdots w^{n}$ on $U_p$ satisfy $w^{n}(q)=0$ for any $q \in U_p
\bigcap K$. In fact, this easily follows by using $w^{n}(exp_{
t\nabla f}(z))=\exp(ht)w^{n}(z)$ to a sequence picked from $W_q
\bigcap U_p$.

As in \cite{B}, one can define a global holomorphic map from $W_p$
to $C$ by extending the local parametrization above. For any point
$q$ distinct from $p$ in $W_p$ we can find a point $q_1$ in $W_p
\bigcap U_p$ such that $q=\exp_{t_1 \nabla f}(q_1)$ for some $t_1$.
Define
$$z^{i}(q)=w^{i}(q_1)  $$ for $1 \leqslant i \leqslant n-1$ and
$$z^{n}(q)=\exp(ht_1)w^{n}(q_1).$$
One can easily check that this definition dose not depend on the
choice of $q_1$. It can also be checked that the holomorphic map $z$
gives a biholomorphic map from $W_p$ to $C$.

The holomorphic line bundle structure of M can be derived from the
above global holomorphic parametrization of $W_p$, thus the theorem
is proved.

\vskip 2 mm

\noindent Before ending this note, we add two remarks here.

\vskip  1 mm

\noindent \textbf{Remark 3.1}. Due to Bryant \cite{B}, the singular
locus of the holomorphic vector field associated to a gradient
K\"{a}hler-Ricci soliton is a disjoint union of nonsingular complex
manifolds, each of which is totally geodesic. In view of this the
assumption of Theorem 3.1 is natural in some sense.

\noindent \textbf{Remark 3.2}. From Theorem 1.1, we can construct
steady gradient K\"{a}hler-Ricci solitons on the canonical line
bundle over a compact K\"{a}hler-Einstein manifold which satisfy all
assumptions in the above theorem Theorem 3.1. It will be interesting
to investigate whether this is the only example.

\end{document}